\documentclass[12pt]{amsart}
\usepackage{amssymb,amsmath,amstext, amsthm,graphicx}
\usepackage{amsmath,amssymb,amscd,amsmath,amstext,mathrsfs,amsfonts,pb-diagram,algorithmic,setspace}%

\newenvironment{Macaulay2}{ \begin{spacing}{0.4} 
\smallskip } { \smallskip 
\end{spacing} }

\newtheorem{theorem}{Theorem}
\theoremstyle{plain}

\newtheorem{algorithm}{Algorithm}

\newtheorem{definition}{Definition}

\newtheorem{problem}{Problem}
\newtheorem{proposition}{Proposition}
\newtheorem{remark}{Remark}

\numberwithin{equation}{section}

\newcommand{\POLYd}{\mathcal{O}_{(d)}}
\newcommand{\ROOTn}{\mathbb{O}^n}
\newcommand{\Pd}{\mathcal{P}_{(d)}}
\newcommand{\CHd}{\mathcal{H}_{(d)}}
\newcommand{\Cn}{\C^n}
\newcommand{\Pn}{\P(\C^{n+1})}

\def\C{{\mathbb C}}
\def\N{{\mathbb N}}

\def\R{{\mathbb R}}

\def\P{{\mathbb P}}
\def\S{{\mathbb S}}

\def\CD{\mathcal{D}}

\newcommand{\mun}{\mu_{\rm norm}}

\newcommand{\Katsura}{\mbox{Katsura}}
\newcommand{\Random}{\mbox{Random}}
\def \LinearHomotopyConstant{0.04804448}

\let\oldmarginpar\marginpar
\renewcommand\marginpar[1]{\-\oldmarginpar[\raggedleft\footnotesize #1]%
{\raggedright\footnotesize #1}}

\title[Certified numerical homotopy tracking]{Certified numerical homotopy tracking}
\author{Carlos Beltr\'an}
        \thanks{C. Beltr\'an. Departamento de Matem\'aticas, Estad\'istica y Computaci\'on, Universidad de Cantabria, Spain
               ({\tt carlos.beltran@unican.es}). Partially supported by MTM2007-62799 and MTM2010-16051, Spanish Ministry of Science (MICINN).}
\author{Anton Leykin}
        \thanks{Anton Leykin. School of Mathematics, Georgia Tech, Atlanta GA, USA ({\tt leykin@math.gatech.edu}).
         Partially supported by NSF grant DMS-0914802}

\date{\today}

\begin{document}
\begin{abstract}
Given a homotopy connecting two polynomial systems we provide a rigorous algorithm for tracking a regular homotopy path connecting an approximate zero of the start system to an approximate zero of the target system. Our method uses recent results on the complexity of homotopy continuation rooted in the alpha theory of Smale. Experimental results obtained with the implementation in the numerical algebraic geometry package of Macaulay2 demonstrate the practicality of the algorithm. In particular, we confirm the theoretical results for random linear homotopies and illustrate the plausibility of a conjecture by Shub and Smale on a good initial pair.
\end{abstract}
\maketitle
The numerical homotopy continuation methods are the backbone of the area of numerical algebraic geometry; while this area has a rigorous theoretic base, its existing software relies on heuristics to perform homotopy tracking. This paper has two main goals:
\begin{itemize}
\item On one hand, we intend to overview some recent developments in the analysis of complexity of polynomial homotopy continuation methods with the view towards a practical implementation. In the last years, there has been much progress in the understanding of this problem. We hereby summarize the main results obtained, writting them in a unified and accesible way.
\item On the other hand, we present for the first time an implementation of a {\em certified} homotopy method which does not rely on heuristic considerations. Experiments with this algorithm are also presented, providing for the first time a tool to study deep conjectures on the complexity of homotopy methods (as Shub \& Smale's conjecture discussed below) and illustrating known --yet somehow surprising -- features about these methods, as equiprobability of the output in the case of random linear homotopy and the average polynomial or quasi--polynomial time of the algorithms studied by several authors.
\end{itemize}
  Our project constructs a certified homotopy tracking algorithm and delivers the first practical implementation of a rigorous path-following procedure. In particular, the case of a {\em linear homotopy} is addressed in full detail in Algorithm~\ref{alg: certified homotopy} of Section~\ref{sec:linear}.

The structure of the paper is as follows. In Section \ref{sec:1} we introduce the basic notations and the general problems that we address. In Section \ref{sec:2} we recall the definition of approximate zero, condition number, and Newton's method, and equip the space of polynomial systems with a Hermitian product. In Section \ref{sec: homotopy method} the main problem solved in this paper is formulated; we describe a certified algorithm to follow a homotopy path. An overview of approaches to finding all the roots of a system is presented in Section \ref{sec:allroots}. In Section \ref{sec: random} we give an algorithm to construct a random linear homotopy with good average complexity. In Section~\ref{sec:implementation} we discuss the implementation of our algorithm. Section~\ref{sec:experiments} demonstrates the practicality of computation with the developed algorithm and discusses experimental data that could be used to obtain intuition, in particular, with regards to a longstanding conjecture of Shub and Smale.

\medskip \noindent {\bf Acknowledgements.} The authors would like to thank Mike Shub for insightful comments; the second author is grateful to Jan Verschelde for early discussions of practical certification issues. This work was partially done while the authors were attending a workshop on the Complexity of Numerical Computation as part of the FoCM Thematic program hosted by the Fields Institute, Toronto. We thank that institution for their kind support.

\section{Description of the problem}\label{sec:1}
Let $n\geq1$. For a positive integer $l\geq1$, let $\mathcal{P}_l=\C_l[X_1,\ldots,X_n]$ be the vector space of all polynomials of degree at most $l$ with complex coefficients and unknowns $X_1,\ldots,X_n$. Then, for a list of degrees $(d)=(d_1,\ldots,d_n)$ let $\Pd=\mathcal{P}_{d_1}\times\cdots\times\mathcal{P}_{d_n}$. Note that elements in $\Pd$ are $n$--tuples $f=(f_1,\ldots,f_n)$ where $f_i$ is a polynomial of degree $d_i$. An element $f\in\Pd$ will be seen both as a vector in some high--dimensional vector space and as a system of $n$ equations with $n$ unknowns.

\begin{problem}\label{prob:affine}
Assuming $f\in\Pd$ has finitely many zeros, find approximately one, several, or all zeros of $f$ in $\C^n$.
\end{problem}

It is helpful to consider the homogeneous version of this problem: For a positive integer $l\geq1$, let $\mathcal{H}_l$ be the vector space of all homogeneous polynomials of degree $l$ with complex coefficients and unknowns $X_0,\ldots,X_n$. Then, for a list of degrees $(d)=(d_1,\ldots,d_n)$ let $\CHd=\mathcal{H}_{d_1}\times\cdots\times\mathcal{H}_{d_n}$. Note that elements in $\CHd$ are $n$--tuples $h=(h_1,\ldots,h_n)$ where $h_i$ is a homogeneous polynomial of degree $d_i$. An element $h\in\CHd$ will be seen both as a vector in some high--dimensional vector space, and as a system of $n$ homogeneous equations with $n+1$ unknowns. Note that if $\zeta\in\C^{n+1}$ is a zero of $h\in\CHd$, then so is $\lambda\zeta$, $\lambda\in\C$. Hence, it makes sense to consider zeros of $h\in\CHd$ as projective points $\zeta\in\Pn$. Abusing the notation, we will denote both a point in $\Pn$ and a representative of the point in $\C^{n+1}$ with the same symbol. Moreover, if necessary, it is implied that the norm of this representative equals 1.

\begin{problem}\label{prob:homogeneous}
Assuming $h\in\CHd$ has finitely many zeros, find approximately one, several or all zeros of $h$ in $\Pn$.
\end{problem}

Let $d=\max\{d_1,\ldots,d_n\}$ and $\CD=d_1\cdots d_n$. Note that $d$ is a small quantity, but in general $\CD$ is an exponential quantity. We denote by $N+1$ the complex dimension of $\CHd$ and $\Pd$ as vector spaces. Namely,
\[
N+1=\sum_{i=1}^n\binom{n+d_i}{d_i}.
\]
There is a correspondence between problems \ref{prob:affine} and \ref{prob:homogeneous}. Given $f=(f_1,\ldots,f_n)\in\Pd$,
\[
f_i=\sum_{\alpha_1+\ldots+\alpha_n\leq d_i}a^i_{\alpha_1,\ldots,\alpha_n}X_1^{\alpha_1}\cdots X_n^{\alpha_n},
\]
we can consider its homogeneous counterpart $h=(h_1,\ldots,h_n)\in\CHd$, where
\[
h_i=\sum_{\alpha_1+\ldots+\alpha_n\leq d_i}a^i_{\alpha_1,\ldots,\alpha_n}X_0^{d_i-(\alpha_1+\cdots+\alpha_n)}X_1^{\alpha_1}\cdots X_n^{\alpha_n},
\]
If $(\eta_1,\ldots,\eta_n)$ is a zero of $f$, then $(1,\eta_1,\ldots,\eta_n)$ is a zero $h$.
Conversely, if $(\zeta_0,\ldots,\zeta_n)\in\Pn$ is a zero of $h$ and $\zeta_0\neq0$ then $\left(\frac{\zeta_1}{\zeta_0},\ldots,\frac{\zeta_n}{\zeta_0}\right)$ is a zero of $f$.

\section{Preliminaries}\label{sec:2}
\subsection{Approximate zeros and Newton's method}

In general, it is hard to describe zeros of $f\in\Pd$ or $h\in\CHd$ exactly. One may ask for points which are ``$\varepsilon$--close'' to some zero, but this is not a very stable concept. The concept of an approximate zero of \cite{Sm86} fixes that gap.

Given $f\in\Pd$, consider the Newton operator associated to $f$,
\[
N(f)(x)=x-Df(x)^{-1}f(x),
\]
where $Df(x)$ is the $n\times n$ derivative matrix of $f$ at $x\in\C^n$, also often called the Jacobian (matrix). Note that $N(f)(x)$ is defined as far as $Df(x)$ is an invertible matrix. We will denote
\[
N(f)^l(x)=N(f)\overset{l}{\overbrace{\circ\cdots\circ}} N(f)(x)
\]
namely, the result of $l$ iterations of Newton's method starting at $x$.
\begin{definition}\label{def:affinezppzero}
We say that $x\in\C^n$ is an approximate zero of $f\in\Pd$ with associated zero $\eta\in\C^n$ if $N(f)^l(x)$ is defined for all $l\geq0$ and
\[
\|N(f)^l(x)-\eta\|\leq\frac{\|x-\eta\|}{2^{2^{l}-1}},\;\;\;\;l\geq0.
\]
\end{definition}
The homogeneous version of Newton's method \cite{Sh93} is defined as follows. Let $h\in\CHd$ and $z\in\Pn$. Then,
\[
N_\P(h)(z)=z-\left(Dh(z)\mid_{z^\perp}\right)^{-1}h(z),
\]
where $Dh(z)$ is the $n\times (n+1)$ Jacobian matrix of $h$ at $z\in\Pn$, and
\[
Dh(z)\mid_{z^\perp}
\]
is the restriction of the linear operator defined by $Dh(z):\C^{n+1}\rightarrow\C^n$ to the orthogonal complement $z^\perp$ of $z$. Hence, $\left(Dh(z)\mid_{z^\perp}\right)^{-1}$ is a linear operator from $\C^n$ to $z^\perp$, and $N_\P(h)(z)$ is defined as far as this operator is invertible. The reader may check that $N_\P(h)(\lambda z)=\lambda N_\P(h)(z)$, namely $N_\P(h)$ is a well--defined projective operator. Note that $N_\P(h)$ may be written in a matrix form
\[
N_\P(h)(z)=z-\binom{Dh(z)}{z^*}^{-1}\binom{h(z)}{0},
\]
which is more comfortable for computations. As before, we denote by $N_\P(h)^l(z)$ the result of $l$ consecutive applications of $N_\P(h)$  with the initial point $z$.
\begin{definition}\label{def:projzppzero}
We say that $z\in\Pn$ is an approximate zero of $h\in\CHd$ with associated zero $\zeta\in\Pn$ if $N_\P(h)^l(z)$ is defined for all $l\geq0$ and
\[
d_R(N_\P(h)^l(z),\zeta)\leq\frac{d_R(z,\zeta)}{2^{2^{l}-1}},\;\;\;\;l\geq0,
\]
\end{definition}
Here $d_R$ is the Riemann distance in $\Pn$, namely
\[
d_R(z,z')=\arccos\frac{|\langle z,z'\rangle|}{\|z\|\,\|z'\|}\in[0,\pi/2],
\]
where $\langle \cdot,\cdot\rangle$ and $\|\cdot\|$ are the usual Hermitian product and norm in $\C^{n+1}$. Note that $d_R(z,z')=d_R(\lambda z,\lambda'z')$ for $\lambda,\lambda'\in\C$, namely $d_R$ is well defined in $\Pn\times\Pn$.
The reader familiar with Riemannian geometry may check that $d_R(z,z')$ is the length of the shortest $C^1$ curve with extremes $z,z'\in\Pn$, when $\Pn$ is endowed with the usual Hermitian structure (see \cite[Page 226]{BlCuShSm98}.)

Let $f\in\Pd$ and let $h\in\CHd$ be the homogeneous counterpart of $f$. In contrast with the case of exact zeros, it may happen that $z=(z_0,\ldots,z_n)$ is an approximate zero of $h$ but still $\left(\frac{z_1}{z_0},\ldots,\frac{z_n}{z_0}\right)$ is not an approximate zero of $f$. In Proposition \ref{prop:fromprojtoaff} we explain how to fix that gap.

\subsection{The Bombieri-Weyl Hermitian product}\label{sec:BW}
In studying Problems \ref{prob:affine} and \ref{prob:homogeneous}, it is very helpful to introduce some geometric and metric properties in the vector spaces $\Pd$ and $\CHd$. We recall now the unitarily--invariant Hermitian product in $\CHd$, sometimes called Kostlan Hermitian product (\cite{BlCuShSm98}) or Bombieri-Weyl Hermitian product (\cite{BeltranShub2009}). Given two polynomials $v,w\in\mathcal{H}_l$,
\[
v=\sum_{\alpha_0+\ldots+\alpha_n=l}a_{\alpha_0,\ldots,\alpha_n}X_0^{\alpha_0}\cdots X_n^{\alpha_n},
\]
\[
w=\sum_{\alpha_0+\ldots+\alpha_n=l}b_{\alpha_0,\ldots,\alpha_n}X_0^{\alpha_0}\cdots X_n^{\alpha_n},
\]
we consider their (Bombieri-Weyl) product
\[
\langle v,w\rangle=\sum_{\alpha_0+\alpha_1+\ldots+\alpha_n=l}\binom{l}{(\alpha_0,\ldots,\alpha_n)}^{-1} a_{\alpha_0,\ldots,\alpha_n}\overline{b_{\alpha_0,\ldots,\alpha_n}},
\]
where $\overline{\;\cdot\;}\,$ is the complex conjugation and
\[
\binom{l}{(\alpha_0,\ldots,\alpha_n)}=\frac{l!}{\alpha_0!\cdots\alpha_n!}
\]
is the multinomial coefficient.

Then, given two elements $h=(h_1,\ldots,h_n)$ and $h'=(h'_1,\ldots,h'_n)$ of $\CHd$, we define
\[
\langle h,h'\rangle=\langle h_1,h_1'\rangle+\cdots+\langle h_n,h_n'\rangle,\;\;\;\;\|h\|=\langle h,h\rangle^{1/2}.
\]
This Hermitian product defines a real inner product in $\CHd$ as usual,
\[
\langle h,h'\rangle_\R=\operatorname{Re}\left(\langle h,h'\rangle\right).
\]
We also define a Hermitian product and the associated norm in $\Pd$ as follows: Given $f,f'\in\Pd$, let $h,h'\in\CHd$ be the homogeneous counterparts of $f,f'$. Then, define
\[
\langle f,f'\rangle=\langle h,h'\rangle,\;\;\;\;\|f\|=\|h\|.
\]
From now on, we will denote by $\S$ the unit sphere in $\CHd$ for this norm, namely
\[
\S=\{h\in\CHd:\|h\|=1\}.
\]
Note that for solving Problem \ref{prob:homogeneous}, we may restrict our input systems $h\in\CHd$ to $h\in\S$, for zeros of a system of equations do not change if the system is multiplied by a non--zero scalar number.

\subsection{The condition number}\label{sec:mu}
The condition number at $(h,z)\in\CHd\times\Pn$ is defined as follows
\[
\mu(h,z)=\|h\|\,\|(Dh(z)\mid_{\,z^\perp})^{-1} \mbox{Diag}(\|z\|^{d_i-1}d_i^{1/2})\|,
\]
or $\mu(h,z)=\infty$ if $Dh(\zeta)\mid_{\,z^\perp}$ is not invertible. Here, $\|h\|$ is the Bombieri-Weyl norm of $h$ and the second norm in the product is the operator norm of that linear operator. Note that $\mu(h,z)$ is essentially equal to the operator norm of the inverse of the Jacobian  $Dh(\zeta)$, restricted to the orthogonal complement of $z$. The rest of the factors in this definition are normalizing factors which make results look nicer and allow projective computations. See \cite{ShSm93b} for more details. Sometimes $\mu$ is denoted $\mun$ or $\mu_{\rm proj}$, but we keep the simplest notation here.

The two following results are versions of Smale's $\gamma$--theorem, and follow from the study of the condition number in \cite{ShSm93b,Shub2007}.

\begin{proposition}\label{prop:aptaffine}\cite[Prop. 4.1]{BePa07}
Let $f\in\Pd$ and let $h\in\CHd$ be its homogeneous counterpart. Let $\eta=(\eta_1,\ldots,\eta_n)\in\C^n$ be a zero of $f$, and let $\zeta=(1,\eta_1,\ldots,\eta_n)\in\Pn$ be the associated zero of $h$. Let $x\in\C^n$ satisfy
\[
\|x-\eta\|\leq\frac{3-\sqrt{7}}{d^{3/2}\mu(h,\zeta)}.
\]
Then, $x$ is an affine approximate zero of $f$, with associated zero $\eta$.
\end{proposition}
\begin{proposition}\label{prop:aptproj}\cite{BeltranTa}
Let $\zeta\in\Pn$ be a zero of $h\in\CHd$ and let $z\in\Pn$ be such that
\[
d_R(z,\zeta)\leq \frac{u_0}{d^{3/2}\mu(h,\zeta)},\;\;\;\;\text{ where $u_0=0.17586$}.
\]
Then $z$ is an approximate zero of $h$ with associated zero $\zeta$.
\end{proposition}
The following result gives a tool to obtain affine approximate zeros from projective ones:
\begin{proposition}\label{prop:fromprojtoaff}\cite[Prop. 4.5]{BePa07}
Let $f\in\Pd$ and let $h\in\CHd$ be its homogeneous counterpart. Let $\eta=(\eta_1,\ldots,\eta_n)\in\C^n$ be a zero of $f$, and let $\zeta=(1,\eta_1,\ldots,\eta_n)\in\Pn$ be the associated zero of $h$. Let $z=(z_0,\ldots,z_n)\in\Pn$ be a projective approximate zero of $h$ with associated zero $\zeta$, such that
\[
d_R(z,\zeta)\leq\arctan\left(\frac{3-\sqrt{7}}{d^{3/2}\mu(h,\zeta)}\right).
\]
($d_R(z,\zeta)\leq \frac{u_0}{d^{3/2}\mu(h,\zeta)}$ suffices.)

Let $z^l=\N_\P(h)^{l}(z)$, where $l\in\N$ is such that
\[
l\geq\log_2\log_2(4(1+\|\eta\|^2)).
\]
Let $x^l=\left(\frac{z^l_1}{z^l_0},\ldots,\frac{z^l_n}{z^l_0}\right)$. Then,
\[
\|x^l-\eta\|\leq\frac{3-\sqrt{7}}{d^{3/2}\mu(f,\eta)}.
\]
In particular, $x^l$ is an affine approximate zero of $f$ with associated zero $\eta$ by Proposition \ref{prop:aptaffine}.
\end{proposition}
Thus, if we have a bound on $\|\eta\|$ and a projective approximate zero of $h$ with associated zero the projective solution $\zeta$, we just need to apply projective Newton's operator $N_\P(h)$ a few times $\lceil\log_2\log_2(4(1+\|\eta\|^2))\rceil$ to get an affine approximate zero of $f$ with associated zero $\eta$. Here, by $\lceil\lambda\rceil$ we mean the smallest integer number greater than $\lambda$, $\lambda\in\R$. Thus, a solution to Problem \ref{prob:affine} follows from a solution to Problem \ref{prob:homogeneous} and a control on the norm of the affine solutions of $f\in\Pd$. The latter can be done either on per case basis or via a probabilistic argument as in \cite[Cor. 4.9]{BePa07}, where it is proved that for $f$ such that $\|f\|=1$ and $\delta\in(0,1)$, we have $\|\eta\|\leq \CD\sqrt{\pi n}/\delta$ with probability greater than $1-\delta$.

From now on we center our attention in Problem \ref{prob:homogeneous}, and we will assume that all the input systems $h$ have unit norm, namely $h\in\S$.

\section{The homotopy method: A one--root finding algorithm} \label{sec: homotopy method}
Let $V=\{(f,\zeta)\in\S\times\Pn:f(\zeta)=0\}$ be the so--called {\em solution variety}. Elements in $V$ are pairs $(system,solution)$. Consider the projection on the first coordinate $\pi:V\rightarrow\S$. The condition number defined above is an upper bound for the norm of the derivative of the local inverse of $\pi$ near $\pi(f,\zeta)$, see for example \cite[Chapter 12]{BlCuShSm98}. In particular, $\pi$ is locally invertible near $(f,\zeta)$ if $\mu(f,\zeta)<\infty$.

Let $t\rightarrow h_t\in\S$, $0\leq t\leq T$ be a $C^1$ curve, and let $\zeta_0$ be a solution of $h_0$. If $\mu(h_0,\zeta_0)<\infty$, then $\pi$ is locally invertible near $h_0$. Thus, there exists some $\varepsilon>0$ such that for $0\leq t<\varepsilon$ the zero $\zeta_0$ can be continued to a zero $\zeta_t$ of $h_t$ in such a way that $t\rightarrow\zeta_t$ is a $C^1$ curve. We call the curve $t\rightarrow(h_t,\zeta_t)$ the lifted curve of $t\rightarrow h_t$. There are two possible scenarios:
\begin{itemize}
\item {\bf Regular:} The whole curve $t\rightarrow h_t$, $0\leq t\leq T$ can be lifted to $t\rightarrow (h_t,\zeta_t)$;
\item {\bf Singular:} There is some $\varepsilon\leq T$ such that $t\rightarrow h_t$ can be lifted for $0\leq t<\varepsilon$, but $\mu(h_t,\zeta_t)\rightarrow\infty$ as $t\rightarrow\varepsilon$.
\end{itemize}

\begin{problem}\label{prob:homotopy}
Create a {\em homotopy continuation algorithm}, a numerical procedure that follows closely the lifted curve. Namely, in the regular case such algorithm's goal is to construct a sequence $0=t_0<t_1<\ldots<t_k=T$ and pairs $(g_i,z_i)\in\S\times\Pn$ such that for all $i=0,\ldots k$ we have  $g_i=h_{t_i}$ and $z_i$ is an approximate zero associated to the zero $\zeta_i$ of $g_i$ with $(g_0,\zeta_0)$ and $(g_i,\zeta_i)$ lying on the same lifted curve.
\end{problem}

The homotopy method that we have in mind would solve the problem above (in the regular case) and would create an infinite sequence $\{t_i\}$ converging to the first singularity on the curve in the singular case.

\begin{remark}\label{rem:homotopy to a singular root}
A homotopy algorithm still may be useful in a singular case where the curve can be lifted for $t\in [0,T)$, which is the scenario, e.g., of a homotopy curve leading to a singular solution. One may use $z_i$ for $t_i$ close to $T$ as an empirical approximation of the singular zero. Approximate zeros (defined before) associated to a singular zero might not exist, since Newton's method loses its quadratic convergence near a singularity.
\end{remark}

Given a $C^1$ curve $t\rightarrow h_t$, we denote $\dot h_t=\frac{d}{dt}h_t$. Namely, $\dot h_t$ is the tangent vector to the curve at $t$. Note that $\dot h_t$ depends on the parametrization of the curve, not only on the geometric object (the arc defined by the curve).

A continuous curve $t\rightarrow h_t\in\S$, $0\leq t\leq T$ is of class $C^{1+Lip}$ if it is of class $C^1$ in $[0,T]$ (i.e. it has a continuous derivative in $(0,T)$ and one--sided derivatives at $t=0$ and $t=T$ making $\dot{h}(t)$ continuous in $[0,T]$), and if the mapping $t\rightarrow \dot h_t$ is a Lipschitz map, namely if there exists a constant $K>0$ such that
\[
\|\dot h_t-\dot h_s\|\leq K|t-s|,\;\;\;\;\;\forall\,t,s\in[0,T].
\]
By Rademacher's Theorem, this implies that the second derivative $\ddot h_t$ exists almost everywhere and is bounded by $\|\dot h_t\|\leq K$.
\subsection{Explicit construction of the homotopy method}\label{sec:ECHM}
A certified homotopy method and its complexity was shown for the first time in \cite{ShSm93b,ShSm94}, at least for the case of linear homotopy. In a recent work \cite{Shub2007}, the theoretical complexity of such methods was greatly improved although no specific algorithm was shown as the choice of the step size was not specified. This last piece can be done in several ways, see \cite{BeltranTa,BurgisserCuckerTa,DedieuMalajovichShub}. We now recall the homotopy method of \cite{BeltranTa}, designed to follow a $C^{1+Lip}$ curve $t\rightarrow h_t \in \S,t\in[0,T]$. We assume that:
\begin{enumerate}
\item We know an approximate zero $z_0$, $\|z_0\|=1$ of $g_0=h_0$, satisfying
\begin{equation}\label{eq:initialcond}
d_R(z_0,\zeta_0)\leq \frac{u_0}{2d^{3/2}\mu(h_0,\zeta_0)},\;\;\;\;\text{ where $u_0=0.17586$},
\end{equation}
for some exact zero $\zeta_0$ of $h_0$.
\item Given $t\in[0,T]$, we can compute $h_t$ and $\dot h_t=\frac{dh_t}{dt}$,
\item We know some real number $H\geq0$ satisfying
\begin{equation}\label{equ:H}
\|\ddot h_t\|\leq d^{3/2}H\|\dot h_t\|^2,
\end{equation}
for almost every $t\in[0,T]$. From now on, we denote
\[
P=\sqrt{2}+\sqrt{4+5H^2}\in\R.
\]
\end{enumerate}
For $i\geq1$, define $(g_{i+1},z_{i+1})$ inductively as follows. Let a representative of $z_i$ be chosen such that $\|z_i\|=1$. Let $s\in[0,T]$ be such that $h_s=g_i$ and let $\dot g_i=\dot h_s\in\CHd$ be the tangent vector to the curve $t\rightarrow h_t$ at $t=s$.
Let
\begin{equation}\label{equ:chi1}
\chi_{i,1}=\left\|\binom{Dg_i(z_i)}{z_i^*}^{-1}\begin{pmatrix}\sqrt{d_1}& & & \\& \ddots & &\\& &\sqrt{d_n}\\&&&1\end{pmatrix}
\right\|,
\end{equation}
\begin{equation}\label{equ:chi2}
\chi_{i,2}=\left(\|\dot g_i\|^2+\left\|\binom{Dg_i(z_i)}{z_i^*}^{-1}\binom{\dot g_i(z_i)}{0}\right\|^2\right)^{1/2},
\end{equation}
and consider
\begin{equation}\label{equ:phi}
\varphi_i=\chi_{i,1}\chi_{i,2}.
\end{equation}
Let
\[
c=\frac{(1-\sqrt{2}u_0/2)^{\sqrt{2}}}{1+\sqrt{2}u_0/2}\left(1-\left(1-\frac{u_0}{\sqrt{2}+2u_0}\right)^{\frac{P}{\sqrt{2}}}\right),
\]
and let $t_i$ be chosen in such a way that
\begin{equation}\label{eq:lowb}
\frac{c}{2Pd^{3/2}\varphi_i}\leq t_i\leq\frac{c}{Pd^{3/2}\varphi_i},
\end{equation}
or $t_i=T-s$ if $\frac{c}{2Pd^{3/2}\varphi_i}\geq T-s$. Note that this last case happens when the step $t_i$ chosen with the formula above takes us beyond the limits of the interval $[0,T]$. The lower bound on (\ref{eq:lowb}) is used to guarantee that the homotopy step is not too small (and thus the total number of steps is not too big!).

Note that in order to compute $\varphi_i$ we must compute the norm of a vector (for $\chi_{i,2}$) and the norm of a matrix (for $\chi_{i,1}$).
However, we only need to do these tasks approximately, for we just need to compute a number in $[\varphi_i,2\varphi_i]$.


In Section \ref{sec:linear} below we describe the value of the constants to be taken in the case of linear homotopy.

 Let $g_{i+1}=h_{s+t_i}$ and let $$z_{i+1}=\frac{N_\P(g_{i+1})(z_i)}{\|N_\P(g_{i+1})(z_i)\|}.$$

This way we generate $(g_1,z_1)$, $(g_2,z_2)$, etc. We stop at $k$ such that $g_k=h_T$, and we output $z_k\in\Pn$.

\subsection{Convergence and complexity of the homotopy method}
The homotopy method is guaranteed to produce an approximate zero of the target system $h=h_T$ if we are in the regular scenario. Moreover, its complexity (number of projective Newton's method steps) is also well understood and attains the theoretical result of \cite{Shub2007}. With the notations above, let
\[
\mathcal{C}_0=\int_0^T\mu(h_t,\zeta_t)\|(\dot h_t,\dot \zeta_t)\|\,dt.
\]
The reader may observe that $\mathcal{C}_0$ is the length of the path $(h_t,\zeta_t)$ in the condition metric, that is the metric in the solution variety $V$ obtained by pointwise multiplying the usual metric (inherited from that of the product $\S\times\Pn$) by the condition number $\mu$.

\begin{theorem}\label{th:stepsize}\cite{BeltranTa}
With the notations and hypotheses above, assume that
\[
d_R(z_0,\zeta_0)\leq\frac{u_0}{2d^{3/2}\mu(h_0,\zeta_0)},\;\;\;u_0=0.17586.
\]
Then, for every $i\geq0$, $z_i$ is an approximate zero of $g_i$, with associated zero $\zeta_i$, the unique zero of $g_i$ that lies in the lifted path $(h_t,\zeta_t)$. Moreover,
\[
d_R(z_i,\zeta_i)\leq\frac{u_0}{2d^{3/2}\mu(h_i,\zeta_i)},\;\;\;\;i\geq1.
\]
If $\mathcal{C}_0<\infty$, there exists $k\geq0$ such that $h_T=g_{k}$. Namely the number of homotopy steps is at most $k$. Moreover,
\[
k\leq \lceil Cd^{3/2}\mathcal{C}_0\rceil,
\]
where
\[
C=\frac{2P}{(1-\sqrt{2}u_0/2 )^{1+\sqrt{2}}} \left(\frac{1}{c}+\frac{1+\sqrt{2}u_0/2 }{\left(1-\sqrt{2}u_0/2 \right)^{\sqrt{2}}}\right).
\]
In particular, if $\mathcal{C}_0<\infty$ the algorithm finishes and outputs $z_k$, an approximate zero of $f=g_k$ with associated zero $\zeta_k$, the unique zero of $f$ that lies in the lifted path $(h_t,\zeta_t)$.
\end{theorem}
\begin{remark}
As $\lceil \lambda\rceil\leq \lambda+1$ for $\lambda\in\R$, we have that the number of steps is at most
\[
1+Cd^{3/2}\mathcal{C}_0.
\]
\end{remark}
\begin{remark}
If the curve $t\rightarrow h_t$ is piecewise $C^{1+Lip}$ we may divide the curve in $L$ pieces, each of them of class $C^{1+Lip}$ and satisfying a.e. $\|\ddot h_t\|\leq d^{3/2}H\|\dot h_t\|^2$ for a suitable $H\geq0$. The algorithm may then be applied to each of these pieces and an upper bound on the total number of steps is at most
\[
L+Cd^{3/2}\mathcal{C}_0.
\]
\end{remark}
\begin{remark}
If more than one approximate zero of $g=f_0$ is known, the algorithm described above may be used to follow each of the homotopy paths starting at those zeros. From Theorem \ref{th:stepsize}, if the approximate zeros of $g$ correspond to different exact zeros of $g$, and if $\mathcal{C}_0$ is finite for all the paths (i.e. if the algorithm finishes for every initial input), then the exact zeros associated with the output of the algorithm correspond to different exact zeros of $f=h_T$.
\end{remark}

\subsection{Linear homotopy}\label{sec:linear}
In the case of linear homotopy, the arc--length parametrization of the path is
\begin{equation}\label{equ: linear homotopy}
t\rightarrow h_t=g\cos(t)+\frac{f-Re(\langle f,g\rangle) g}{\sqrt{1-Re(\langle f,g\rangle)^2}}\sin(t),\;\;\;t\in\left[0,T\right],
\end{equation}
where
\[
T=\arcsin\sqrt{1-Re\langle f,g\rangle^2}=\operatorname{distance}(g,f)\in[0,\pi].
\]
Note that this is a $C^\infty$ parametrization so in particular it is $C^{1+Lip}$. From \cite[Section 2.2]{BeltranTa}, in this case we may take the following value of $c/P$ in the description of the algorithm,
\[
\frac{c}{P}=\LinearHomotopyConstant...
\]

The procedure of certified tracking for a linear homotopy is presented by Algorithm~\ref{alg: certified homotopy}.

\begin{algorithm} \label{alg: certified homotopy} $z_*=TrackLinearHomotopy(f,g,z_0)$
\begin{algorithmic}[1]
\REQUIRE $f,g\in\S$; $z_0$ is an approximate zero of $g$ satisfying (\ref{eq:initialcond}).
\ENSURE $z_*$ is an approximate zero of $f$ associated to the end of the homotopy path starting at the zero of $g$ associated to $z_0$
and defined by the homotopy (\ref{equ: linear homotopy}).
\STATE $i \leftarrow 0$; $s_i=0$.
\WHILE {$s_i \neq T$}
\STATE Compute
\[
\dot g_i \leftarrow \dot h_s=-g\sin(s)+\frac{f-Re(\langle f,g\rangle) g}{\sqrt{1-Re(\langle f,g\rangle)^2}}\cos(s).
\]
at $s=s_i$.
\STATE Determine $\varphi_i$ using (\ref{equ:chi1}), (\ref{equ:chi2}), and (\ref{equ:phi}).
\STATE Let $t_i$ be any number satisfying
\[
\frac{\LinearHomotopyConstant}{2d^{3/2}\varphi_i}\leq t_i\leq\frac{\LinearHomotopyConstant}{d^{3/2}\varphi_i}.
\]
\IF{$t_i>T-s$}
\STATE $t_i\leftarrow T-s$.
\ENDIF
\STATE $s_{i+1} \leftarrow s_i + t_i $; $g_{i+1} \leftarrow h_{s_{i+1}}$; $z_{i+1}=\frac{N_\P(g_{i+1})(z_{i})}{\|N_\P(g_{i+1})(z_{i})\|}$.
\STATE $i \leftarrow i + 1$.
\ENDWHILE
\STATE $z_* \leftarrow z_T$.
\end{algorithmic}
\end{algorithm}

The bound on the number of steps in Algorithm~\ref{alg: certified homotopy} given by Theorem~\ref{th:stepsize} is
\begin{equation}\label{eq:k}
k\leq \lceil 71d^{3/2}\mathcal{C}_0\rceil.
\end{equation}

\section{Finding all roots}\label{sec:allroots}

Let us consider polynomial functions in $\POLYd$, where $\POLYd$ is one of $\{\Pd, \CHd, \S\}$ with zeros in $\ROOTn$, where $\ROOTn$ is either $\Cn$ or $\Pn$.

Consider a homotopy $t\rightarrow h_t \in \POLYd$, $t\in[0,T]$, connecting the {\em target system} $h_T$ and the {\em start system} $h_0$ along with a set of {\em start solutions} $Z_0 \subset h_0^{-1}(0) \subset \ROOTn$.

Suppose the homotopy curve $t\to h_t$ can be lifted to $t\to (h_t,\zeta_t)\in \POLYd\times\ROOTn$, $t\in[0,T]$ such that the projection map $\pi: \POLYd\times\ROOTn \to \POLYd$ is locally invertible at any $t\in [0,T)$. A {\em homotopy path} is defined as the projection of such lifted curve onto the second coordinate. If the map $\pi$ is locally invertible at $t=T$ as well, then the path is called {\em regular}.

The homotopy $t\rightarrow h_t$ is called {\em optimal} if every $\zeta_0\in Z_0$ is the beginning of a regular homotopy path.
If every solution of $h_T$ is the (other) end of the homotopy path beginning at some $\zeta_0\in Z_0$ then we call the homotopy {\em total}.

The homotopy method described in Section \ref{sec: homotopy method} can also be applied to follow all homotopy paths of an optimal homotopy in $\S$. Namely, if we have a $C^{1+Lip}$ curve $t\rightarrow h_t$ satisfying the hypotheses of Theorem \ref{th:stepsize} and start with approximate zeros $\tilde Z_0$ associated to a set of solutions $Z_0 \subset h_0^{-1}(0) \subset \Pn$, then we may perform the algorithm for each of the initial pairs $(h_0,z_0)$, $z_0\in Z_0$. By Theorem \ref{th:stepsize}, the output will be approximate zeros associated to $\sharp (Z_0)$ distinct solutions of $h_T$.

The area of {\em numerical algebraic geometry} (see, e.g., \cite{Sommese-Wampler-book-05}) relies on the ability to reliably track optimal homotopies and find {\em all} roots of a given 0-dimensional polynomial system of equations in $\POLYd$. To accomplish that one has to arrange a total homotopy.

\subsection{Total-degree homotopy}\label{sec:totaldegree}
For a target system in $h_T\in\Pd$, $(d)=(d_1,\ldots,d_n)$, define a {\em total degree linear homotopy} to be
\begin{equation}\label{equ:totaldegree}
t\to h_t = (T-t) h_0 + \gamma t h_T,\ \gamma \in \C^*,\ t\in[0,T],
\end{equation}
where the start system is
\begin{equation}\label{equ:start-totaldegree}
  h_0 = (x_1^{d_1} - 1, \ldots, x_n^{d_n} - 1) \in \Pd.
\end{equation}
There are {\em total degree} many, i.e., $d_1\cdots d_n$, zeros of $h_0$ that one can readily write down.
\begin{proposition}\label{prop:Bezout}
Assume that $h_T$ has a finite number of zeros, and let $Z_0$ be the set of zeros of $h_0$ above. Then for all but finitely many values of the constant $\gamma$ the homotopy (\ref{equ:totaldegree}) is total.

If the target system $h_T\in\Pd$ has total-degree many solutions, then (for a generic $\gamma$) the homotopy is optimal.
\end{proposition}

If the target system $h_T\in\Pd$ has fewer than total degree many solutions then:
 \begin{itemize}
   \item some solutions of the target system may be multiple (singular);
   \item in case of $\ROOTn = \Cn$, some of the homotopy paths may diverge (to infinity) when approaching $t=T$.
 \end{itemize}

\begin{remark}\label{rem:regularization}
To compute singular solutions one may track regular homotopy paths to $t=T-\varepsilon$ for a small $\varepsilon>0$ (as in Remark \ref{rem:homotopy to a singular root}) and then use either {\em singular endgames} \cite[Section 10.3]{Sommese-Wampler-book-05} or {\em deflation} \cite{LVZ,LVZ-higher}.

To avoid diverging paths one may homogenize the homotopy passing from $\Pd$ to $\CHd$.
\end{remark}

\begin{remark}\label{rem: thm application}
Normalizing an optimal homotopy $t\to h_t$ with respect to the Bombieri-Weyl norm brings the system from $\CHd$ to $\S$. In case of a linear homotopy, arc-length re-parametrization enables an application of Algorithm~\ref{alg: certified homotopy} for each start solution.
\end{remark}

\subsection{Other homotopy methods}
There are other ways to obtain all solutions with homotopy continuation that exploit either sparseness or special structure of a given polynomial system, here we list a few:
\begin{itemize}
\item Polyhedral homotopy continuation based on \cite{HuberSturmfels:PolyhedralHomotopies} allows to recover all solutions of a sparse polynomial system in the torus $(\C^*)^n$.
\item In many cases presented with a parametric family of polynomial systems it is enough to solve one system given by a generic chioice of parameters. Then, given another system in the family, the chosen generic system may be used as a start system in the so-called {\em coefficient-parameter} or {\em cheater's} homotopy \cite[Chapter 7]{Sommese-Wampler-book-05} recovering all solutions of the latter.
\item Special homotopies: e.g., {\em Pieri homotopies} coming up in Schubert calculus \cite{HuberSottileSturmfels:NumericalSchubertCalculus} are total and optimal by design.
\end{itemize}
For the purpose of this paper we assume that some regularization procedure (see Remark \ref{rem:regularization}) has been applied to make these homotopies optimal and they are brought to $\S$ as in Remark \ref{rem: thm application} .

\section{Random linear homotopy and polynomial time}\label{sec: random}

Suppose, given a regular system $f\in \CHd$, we would like to construct an initial pair $(g,\zeta_0)$ in a random fashion so that every root of $f$ is equally likely to be at the end of the linear homotopy path determined by this initial pair. A simple solution to this problem would be to take $g$ to be the start system (\ref{equ:start-totaldegree}) of the total-degree homotopy -- or its homogenized version --  and pick $\zeta_0$ from the start solutions with uniform probability distribution on the latter. It has been very recently proved by B\"urgisser and Cucker \cite{BurgisserCuckerTa} that this is a pretty good candidate for the linear homotopy starting pair, as the total average number of steps for each path is $O(d^3Nn^{d+1})$ that is $O(N^{\log(\log(N))})$, hence close to polynomial in the input size, mainly when $n>>d$.

\smallskip

In \cite{BePa:FoCM,BePa07,BeltranPardoFLH} a probabilistic way to choose the initial pair was proposed. We now center our attention in the most recent of these works \cite{BeltranPardoFLH} where it is proved that, if the initial pair $(g,\zeta_0)$ is chosen at random (with a certain probability distribution), then the average number of steps performed by the algorithm described in Section \ref{sec: homotopy method} is $O(d^{3/2}nN)$, thus almost linear in the size of the input. It is also proved that in this way we obtain an approximation of a zero of $f$, so that each of the zeros of $f$ are equiprobable if $f$ has no singular solution. In \cite{BeltranShub2009} it is seen that some higher moments (in particular, the variance) of that algorithm are also polynomial in the size of the input. In this section we describe in detail how the process of randomly choosing $(g,\zeta_0)$ works and we recall the main results of \cite{BeltranPardoFLH,BeltranShub2009}.

Given $\zeta\in\Pn$, we consider the set
\[
R_\zeta=\{\tilde{h}\in\CHd:\tilde{h}(\zeta)=0,D\tilde{h}(\zeta)=0\}.
\]
Note that $R_\zeta$ is defined as the set of polynomials in $\CHd$ whose coefficients (in the usual monomial basis) satisfy $n^2+2n$ linear homogeneous equalities. Thus, $R_\zeta$ is a vector subspace of $\CHd$. Moreover, let $e_0=(1,0,\ldots,0)^T$. Then, $R_{e_0}$ is the set of polynomial systems  $\tilde{h}=(\tilde{h}_1,\ldots,\tilde{h}_n)\in\CHd$ such that all the coefficients of $\tilde{h}_i$ containing $X_0^{d_i}$ or $X_0^{d_i-1}$ are zero, namely
\[
\tilde{h}_i=X_0^{d_i-2}p_{2,i}(X_1,\ldots,X_n)+X_0^{d_i-3}p_{3,i}(X_1,\ldots,X_n)+\cdots,
\]
where $p_{k,i}$ is a homogeneous polynomial of degree $k$ with unknowns $X_1,\ldots,X_n$.

Recall that $N+1$ is the (complex) dimension of $\CHd$. The process of choosing $(g,\zeta_0)$ at random is as follows:
\begin{enumerate}
\item Let $(M,l)\in\C^{n^2+n}\times\C^{N+1-n^2-n}=\C^{N+1}$ be chosen at random with the uniform distribution in
\[
B(\C^{N+1})=\{r\in\C^{N+1}:\|r\|_2\leq1\},
\]
where $\|\cdot\|_2$ is the usual Euclidean norm in $\C^{N+1}$. Thus, $M$ is a $(n^2+n)$--dimensional complex vector, that we consider as a $n\times (n+1)$ complex matrix. Note that choosing $\|(M,l)\|_2\leq1$ implies that $\|M\|_F\leq1$ and indeed the expected value of $\|M\|_F^2$ is $\frac{n^2+n}{N+2}$. At this point we can discard $l$ and just keep $M$. Note that this procedure is different from just choosing a random matrix, as it induces a certain distribution in the norm of the matrix which is precisely the one that we are interested in. Hence, choosing $(M,l)$ in the unit ball and then discarding $l$ is not a fool job!
\item With probability $1$, the choice above has produced a matrix $M$ whose kernel has complex dimension $1$. Let $\zeta_0$ be a unit norm element of $Ker(M)$, randomly chosen in $Ker(M)$ with the uniform distribution (we may just obtain any such $\zeta_0$ by solving $M\zeta_0=0$ with our preferred method, and then multiply $\zeta_0$ by a uniformly chosen random complex number of modulus $1$). Let $V$ be any unitary matrix such that $V^*\zeta_0=e_0$. Choose a system $\tilde{h}$ at random in the unit ball (for the Bombieri--Weyl norm) of $R_{e_0}$.  Then, consider $h=\tilde{h}\circ V^*$. (This last procedure is equivalent to choosing a system at random with the uniform distribution in $B(R_{\zeta_0})=\{h\in R_{\zeta_0}:\|h\|\leq1\}$.)
\item Let $\hat{g}\in\CHd$ be the polynomial system defined by
\[
\hat{g}(z)=\sqrt{1-\|M\|_F^2} h(z)+\begin{pmatrix} \langle z,\zeta_0\rangle^{d_1-1}\sqrt{d_1}&&\\&\ddots&\\&&\langle z,\zeta_0\rangle^{d_n-1}\sqrt{d_n}\end{pmatrix}Mz
\]
\item Let
\[
g=\frac{\hat{g}}{\|\hat{g}\|}.
\]
Then, we have chosen $(g,\zeta_0)$ and the reader may check that $g(\zeta_0)=0$, so $\zeta_0$ is an exact zero of $g$.
\end{enumerate}

Consider the randomized algorithm defined as follows:
\begin{enumerate}
\item Input $f\in\S$
\item Choose $(g,\zeta_0)$ at random with the process described above
\item Consider the path
\[
t\rightarrow h_t= g\cos(t)+\frac{f-Re(\langle f,g\rangle) g}{\sqrt{1-Re(\langle f,g\rangle)^2}}\sin(t),\;\;\;t\in\left[0,T\right],
\]
where $T=\arcsin\sqrt{1-Re\langle f,g\rangle^2}$, and note that $h_0=g,h_T=f$. Use Algorithm \ref{alg: certified homotopy} to follow the path $h_t$ and output an approximate zero of $f$.
\end{enumerate}
For given $f\in\S$, let $NS(f)$ be the expected number of homotopy steps performed by this algorithm, on input $f\in\S$. We have seen in (\ref{eq:k}) that
\[
NS(f)\leq \lceil71d^{3/2}\mathcal{C}_0\rceil
\]
The main theorems of \cite{BeltranPardoFLH,BeltranShub2009} are now summarized as follows.
\begin{theorem}\label{th:complexityhomotopy}
If $f\in\S$ is such that every zero of $f$ is non--singular (thus, $f$ has exactly $\mathcal{D}=d_1\cdots d_n$ projective zeros), then:
\begin{itemize}
\item The algorithm above finishes with probability $1$ on the choice of $(g,\zeta_0)$, and
\item Every zero of $f$ is equally probable as the exact zero associated with the output of the algorithm (which is an approximate zero of $f$).
\end{itemize}
Assuming $f\in\S$ is chosen at random with the uniform distribution on $\S$, the expected value and variance of $NS(f)$ satisfy
\[
{\rm E}(NS(f))\leq C_1nNd^{3/2},\;\;\;\;{\rm Var}(NS(f))\leq C_2n^2N^2d^3\ln(\mathcal{D}),
\]
where $C_1$ and $C_2$ are universal constants. One may choose $C_1=71\pi/\sqrt{2}$.
\end{theorem}

Note that this theorem not only gives a uniform distribution of the probability of producing any given root of a regular system, but also gives a good expected running time, with a number of steps which is almost linear in the size of the input.

An algorithm for finding all solutions of a system $f$ with regular zeros follows from Theorem~\ref{th:complexityhomotopy}: repeatedly create and follow random homotopies to find one root of the system until total-degree many roots are found. Assuming a very small probability of finding less than the all the roots, it suffices to choose $\mathcal{D}\log\mathcal{D}$ such random homotopies. Thus, the expected number of steps of the proposed procedure is $O(d^{3/2}nN\mathcal{D}\log\mathcal{D})$, which grows fast as the total degree of the system increases. This fast growth is necessary if we are attempting to find all the $\CD$ solutions of the system. The bound $O(d^{3/2}nN\mathcal{D}\log\mathcal{D})$ is the smallest proven value for the complexity of finding all roots of a system. However, this algorithm may not be the most practical one. Using the na\"ive start system (\ref{equ:start-totaldegree}) should require, according to \cite{BurgisserCuckerTa} an average number of steps $O(d^{3/2}n^{d+1}N\mathcal{D})$ which is a bigger bound that $O(d^{3/2}nN\mathcal{D}\log\mathcal{D})$, but guarantees that just $\CD$ homotopy paths have to be followed.

\section{Implementation of the method}\label{sec:implementation}

The computer algebra system Macaulay2 -- to be more precise, NAG4M2 (internal name: {\tt NumericalAlgebraicGeometry}) package \cite {Leykin:NAG4M2} -- hosts the implementation of Algorithm~\ref{alg: certified homotopy}, which is the first implementation of certified homotopy tracking in a numerical polynomial homotopy continuation software. The current implementation is carried out with standard double floating point arithmetic without analyzing effects of round-off errors. For a variant of the algorithm that facilitates rigorous error control see~\cite{Beltran-Leykin:TRC}.

\subsection{NAG4M2: User manual}
There are several functions that we would like to describe here. First let us give an example of launching {\tt track} procedure with the certified homotopy tracker:
{\small
\begin{Macaulay2}
\begin{verbatim}
i1 : loadPackage "NumericalAlgebraicGeometry";

i2 : R = CC[x,y,z];

i3 : T = {x^2+y^2-z^2, x*y};

i4 : (S,x0) = totalDegreeStartSystem T;

i5 : x1 = first track(S,T,x0,
                      Predictor=>Certified,Normalize=>true)

o5 = {.00000207617, -.706804, .70744}

o5 : Point

i6 : x1.NumberOfSteps

o6 = 129
\end{verbatim}
\end{Macaulay2}
} 
\smallskip
\noindent The values for the optional arguments {\tt Predictor} and {\tt Normalize} specify that the certified homotopy tracking is performed and the polynomial systems are normalized to the unit sphere $\S$.
In this particular example, {\tt totalDegreeStartSystem} creates an initial pair based on the homogenization of the system described in~\ref{equ:start-totaldegree} and {\tt track} follows the linear homotopy starting at this initial pair and finishing at the given target system. 

The user can also get a good initial pair~(\ref{equ: conjecture}) discussed below with the function {\tt goodInitialPair} as well as a random pair of start system and solution as described in Section~\ref{sec: random} with {\tt randomInitialPair}.

It is possible for {\tt track} to return a solution marked as {\em failure}. This happens when the step size becomes smaller than the threshold set by the optional parameter {\tt tStepMin}, which has the default value of $10^{-6}$.

\subsection{Uncertified homotopy continuation}
All existing software, such as HOM4PS2 \cite{HOM4PSwww}, Bertini \cite{Bertini}, and PHCpack \cite{V99}, utilize algorithms based on alternating {\em predictor} and {\em corrector} steps. Here is a summary of operations performed at a point of continuation sequence $t \in [0,T]$ starting with a pair $(h_t,x_t)$ where $x_t$ approximates some zero $\eta_t$ of $h_t$:
\begin{enumerate}
  \item \label{heuristic1} Decide heuristically on the step size $\Delta t$ that predictor should take;
  \item \label{heuristic2} Use a predictor method, i.e., one of the methods for numerical integration of the system of ODEs
  $$
  \dot z = - (Dh_t)^{-1}\, \dot h_t
  $$
  to produce an approximation of $\zeta_{t+\Delta t}$, a solution of $h_{t+\Delta t}$;
  \item \label{heuristic3} Apply the corrector: perform a fixed number $l$ of iterations of Newton's method to obtain a refined approximation $x_{t+\Delta t} = N(h_{t+\Delta t})^l(x_{t+\Delta t})$;
  \item \label{heuristic4} If the estimated error bound in step \ref{heuristic3} is larger than a predefined tolerance, decrease $\Delta t$ and go to step \ref{heuristic1}.
\end{enumerate}
After tuning the parameters, e.g., tolerances values, the application of described heuristics often produces correct solutions. 

\subsection{Certified vs Heuristic}
We can imagine several ``unfortunate'' scenarios when two distinct homotopy paths come too close to each other. Consider sequences $z_0,z_{t_1},\ldots,z_{t_k}$ and $z'_0,z'_{t'_1},\ldots,z'_{t'_{k'}}$ created by an uncertified algorithm in an attempt to approximate these two paths:
\begin{itemize}
  \item If there are subsequences in two sequences that approximate a part of the same path then this is referred to as {\em path jumping}.
  \item {\em Path crossing} happens when the sequences jump from one path to the other, but there is no common path segment that they approximate.
\end{itemize}
While path jumping can be detected, in principle, {\em a posteriori} and the continuation rerun with tighter tolerances and smaller step sizes, the path crossing can not be determined easily.

Path crossing does not result in an incorrect set of target solutions; however, for certain homotopy-based algorithms such as {\em numerical irreducible decomposition}~\cite{SVW1} and applications relying on {\em monodromy} computation such as~\cite{Leykin-Sottile:HoG} the order of the target solutions is crucial. Therefore, one not only needs to certify the end points of homotopy paths, but also has to show that the approximating sequences follow the same path from start to finish. The certification of the sequence produced in Section \ref{sec: homotopy method} provided by Theorem \ref{th:stepsize} gives such guarantee.

In certain cases the target solutions obtained by means of uncertified homotopy continuation can be rigourously certified after all of them are obtained. For instance, suppose a target system $h_T\in \CHd$ has distinct regular solutions in $\Pn$, then there are total degree many of them. Suppose some procedure provides total degree many approximations to solutions. If a bound on $\max\{\mu(h_T,\zeta)\,|\,\zeta\in h_T^{-1}(0)\}$ is known, then using Proposition~\ref{prop:aptproj} these approximations may be certified as distinct numerical zeros, thus certifying that all solutions have been found. If no such bound is known, one may still try to prove that the zeros are different by means of Smale's $\alpha$--theorem \cite{Sm86} (see \cite{Hauenstein-Sottile:alphaCertified}). However, these procedures do not detect if path crossing has occurred.

\section{Experimental results}\label{sec:experiments}

The developed and implemented algorithm provided us with a chance to conduct experiments that illuminated several aspects in the complexity analysis of solving polynomial systems via homotopy continuation.

\subsection{Practicality of certified tracking}
Our experiments in this section were designed to explore how practical the certified tracking provided by Algorithm~\ref{alg: certified homotopy} is. As was already mentioned, the proposed certified procedure makes sense only for a regular homotopy. Moreover, in nearly singular examples the certified homotopy is bound to show bad performance due to steps being minuscule at the end of paths, which is mandated by (\ref{eq:lowb}).

In the table below we give the data produced by tracking of total-degree homotopy that is optimal for the chosen examples:
\begin{itemize}
  \item $\Random_{(d_1,\ldots,d_n)}$: a random system in $\S \subset \CHd$ with uniform distribution;
  \item $\Katsura_n$:  a classical benchmark with one linear and $n-1$ quadratic equations in $n$ variables.
\end{itemize}

For every experiment we provide the number of solutions, the average number of steps per homotopy path both for the certified algorithm (C) and one of the best heuristic procedures (H) implemented in Macaulay2.

\begin{center}
\begin{tabular}{|l||c|c|c|}
\hline
system & \#sol. & \#steps/path (C) & \#steps/path (H)\\
\hline
$\Random_{(2,2)}$ & 4 & 198.5 & 31\\
$\Random_{(2,2,2)}$ & 8 & 370.125 & 23\\
$\Random_{(2,2,2,2)}$ & 16 & 813.812 & 44.375\\
$\Random_{(2,2,2,2,2)}$ & 32 & 1542.5 & 48.5312\\
$\Random_{(2,2,2,2,2,2)}$ & 64 & 2211.58 & 58.5312\\
$\Katsura_3$ & 4 & 569.5 & 25.75\\
$\Katsura_4$ & 8 & 1149.88 & 41.5\\
$\Katsura_5$ & 16 & 1498.38 & 39.0625\\
$\Katsura_6$ & 32 & 2361.81 & 55.5625\\
\hline
\end{tabular}
\end{center}

One step in a heuristic algorithm takes more work than that of the certified tracker: there is a predictor and several corrector steps performed and, if unsuccessful, new step size chosen only to repeat the procedure. Despite that the heuristic approach leads to much smaller computational time for larger systems: it could be concluded from the table above that the number of successful heuristic steps does not grow fast with degree of the system and the number of variables (in comparison to certified tracking).

Of course, it is clear that if we want to run a certified, non heuristic method as the one we propose, we will need more computational time.

\subsection{A conjecture by Shub and Smale}
In \cite{ShSm94}, a pair 
\begin{equation}\label{equ: conjecture}
g(x)=\begin{cases}d_1^{1/2}x_0^{d_1-1}x_1\\ \vdots\\ d_n^{1/2}x_0^{d_n-1}x_n\end{cases},\;\;\;\;e_0=\begin{pmatrix}1\\0\\ \vdots\\0\end{pmatrix}.
\end{equation} 
was conjectured to be a good starting pair for the linear homotopy. More exactly, let
\[
E_{good}=E(\sharp(steps)\text{ to solve f with lin. homotopy starting at }(g,e_0)),
\]
where the expectation is taken for random $f\in\S$. Then, the conjecture in \cite{ShSm94} can be writen as\footnote{The original pair suggested by Shub and Smale had no $d_i^{1/2}$ factors as the one here. As done in other papers, we add those factors here to optimize the condition number $\mu(g,e_0)$.}
\begin{equation}\label{eq:SS94}
E_{good}\leq\text{ a small quantity, polynomial in $N$},
\end{equation}

The following experimental data was obtained by running a linear homotopy connecting the pair $(g,e_0)$ as in (\ref{equ: conjecture}) to a random system in $\S\subset\CHd$ with $d_i=2$ for $i=1,\ldots,n$. We compare the values to that of $B(n,d,N) = 71\pi d^{3/2}nN/\sqrt{2}$ which according to Theorem \ref{th:complexityhomotopy} is a bound for the average number of steps.

{\small
$$
\begin{array}{|c||c|c|c|c|c|c||c|}
\hline
 n &    E_{good}  &  Var_{good}  &   E_{total}   &  Var_{total} &  E_{rand}   &    Var_{rand}  &   B(n,d,N)        \\
\hline
4  &    962.051   &3.2\cdot 10^5&  1263.72      &4.3\cdot 10^5&  1622.29    &  6.8\cdot 10^5    & 1.0\cdot 10^5 \\
5  &    1524.6    &6.9\cdot 10^5&  2130.54      &1.2\cdot 10^6&  2728.3     &  1.7\cdot 10^6    & 2.3\cdot 10^5 \\
6  &    2258.33   &1.3\cdot 10^6&  3129.56      &2.2\cdot 10^6&  4137.16    &  3.5\cdot 10^6    & 4.5\cdot 10^5 \\
7  &    3130.83   &2.3\cdot 10^6&  4530.55      &4.5\cdot 10^6&  5743.32    &  5.5\cdot 10^6    & 7.8\cdot 10^5 \\
8  &    4154.38   &3.9\cdot 10^6&  5967.57      &6.7\cdot 10^6&  8048.94    &  1.0\cdot 10^7    & 1.2\cdot 10^6 \\
9  &    5488.93   &7.0\cdot 10^6&  8013.71      &1.1\cdot 10^7&  10482.1    &  1.6\cdot 10^7    & 1.9\cdot 10^6 \\
10 &    6871.35   &1.0\cdot 10^7&  10071        &1.4\cdot 10^7&  13477.5    &  2.2\cdot 10^7    & 2.9\cdot 10^6 \\
11 &    8622      &1.2\cdot 10^7&  12996.1      &2.8\cdot 10^7&  17193.3    &  3.5\cdot 10^7    & 4.2\cdot 10^6 \\
12 &    10413.3   &2.0\cdot 10^7&  15115.4      &2.8\cdot 10^7&  20761.3    &  4.6\cdot 10^7    & 5.8\cdot 10^6 \\
13 &    12447.1   &2.6\cdot 10^7&  18744.5      &4.3\cdot 10^7&  25646.5    &  6.3\cdot 10^7    & 7.9\cdot 10^6 \\
14 &    14769.9   &3.3\cdot 10^7&  22317.1      &6.1\cdot 10^7&  29596.7    &  9.1\cdot 10^7    & 1.0\cdot 10^7 \\
15 &    17255.7   &4.4\cdot 10^7&  26017.7      &7.3\cdot 10^7&  35582.6    &  1.2\cdot 10^8    & 1.4\cdot 10^7 \\
16 &    20959.7   &5.9\cdot 10^7&  30063.9      &1.0\cdot 10^8&  42098.9    &  1.5\cdot 10^8    & 1.7\cdot 10^7 \\
17 &    23589.4   &7.5\cdot 10^7&  35403.1      &1.3\cdot 10^8&  48024.5    &  1.7\cdot 10^8    & 2.2\cdot 10^7 \\
18 &    27400.9   &9.6\cdot 10^7&  40242.5      &1.5\cdot 10^8&  54955.4    &  2.3\cdot 10^8    & 2.7\cdot 10^7 \\
19 &    29930.3   &1.0\cdot 10^8&  46502.2      &2.3\cdot 10^8&  62855.2    &  2.9\cdot 10^8    & 3.4\cdot 10^7 \\
20 &    34374.2   &1.4\cdot 10^8&  51730.2      &2.3\cdot 10^8&  71242.5    &  3.5\cdot 10^8    & 4.1\cdot 10^7 \\
\hline
\end{array}
$$
}

For each value of $n$ we have generated 1000 random systems in $\S$ with a uniform probability distribution.
The values $E_{good}$ and $Var_{good}$ are estimated expected value and variance of the number of steps taken by Algorithm~\ref{alg: certified homotopy} for the initial pair in (\ref{equ: conjecture}); $E_{rand}$ and $Var_{rand}$ refer to those for the {\em random} initial pair; $E_{total}$ and $Var_{total}$ refer to those for the homogeneous version of the total--degree homotopy system of Section \ref{sec:totaldegree} containing all the roots of unity (the choice of the root is irrelevant for symmetry reasons). Namely, the pair
\begin{equation}\label{eq:totdegreehom}
h_0=(X_1^{d_1}-X_0^{d_1},\ldots,X_n^{d_n}-X_0^{d_n}),\;\;\; \zeta_0=(1,\ldots,1).
\end{equation}

\begin{figure}[b]
  \begin{center}
    \includegraphics[width=14cm, keepaspectratio]{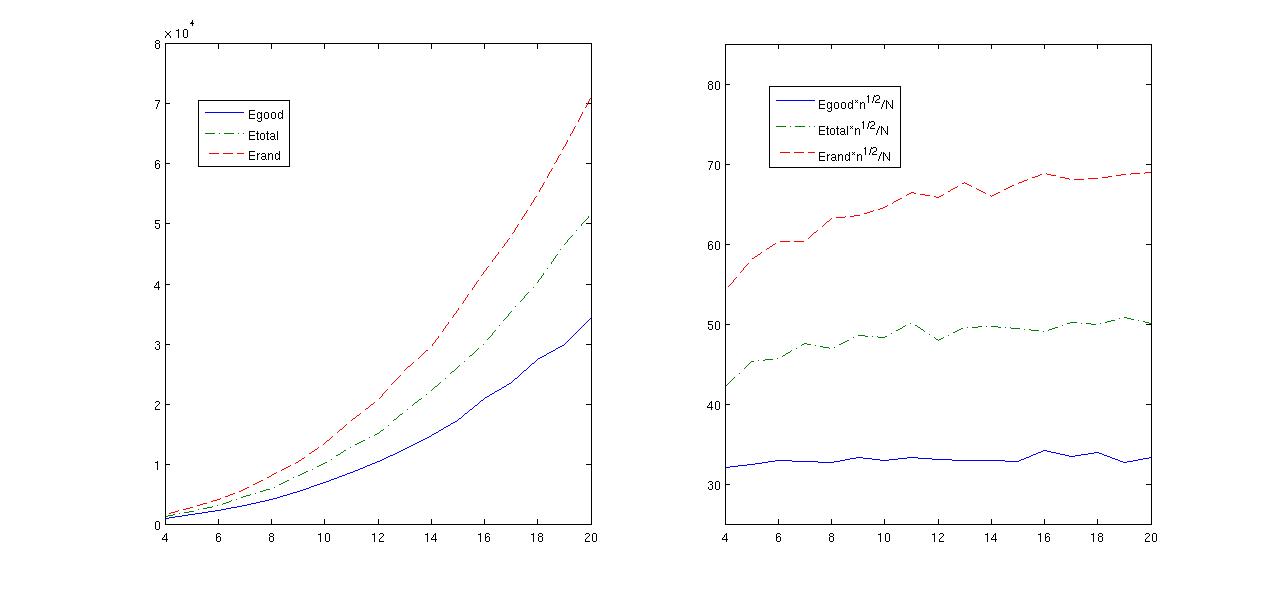}
  \end{center}
  \caption{In the first figure, we have plotted the experimental values obtained for $E_{good},E_{rand}$ and $E_{total}$ for $n=4,\ldots,20$. In the second one we plot $\frac{E_{good}n^{1/2}}{N}$, $\frac{E_{total}n^{1/2}}{N}$ and $\frac{E_{rand}n^{1/2}}{N}$ for $n=4,\ldots,20$}
  \label{figure:1}
\end{figure}

The table above and Figure~\ref{figure:1} below suggest two conclusions for the case of degree two polynomials:
\begin{itemize}
\item The random homotopy seems to take approximately double number of steps than the homotopy with initial pair (\ref{equ: conjecture}). The total degree homotopy lies in-between.
\item The average number of steps in the three cases seem to grow as $Constant \cdot\frac{N}{\sqrt{n}}$ with $Constant\approx 35,50,70$ for $E_{good}$, $E_{total}$ and $E_{rand}$ respectively.
\end{itemize}
This experiment thus confirms the conjecture by Shub and Smale and moreover it suggests a more specific form, suggesting that the same bound given for random homotopy should hold for the conjectured pair:
\begin{equation}\label{eq:SS94B}
E_{good}\leq CnNd^{3/2},
\end{equation}
with $C$ a constant. 
We also extend this conjecture to the case of the initial pair total--degree homotopy pair $(h_0,\zeta_0)$ of Equation (\ref{eq:totdegreehom}) above:
\[
E_{total}\leq CnNd^{3/2}.
\]
Moreover, as pointed out above, in the case of degree $2$ systems, we obtained experimentally a better bound
\[
E_{good},E_{total},E_{rand}\leq\frac{CN}{\sqrt{n}},
\]
$C$ a constant. The difference between the experimentally observed value and the theoretical bound in the case of randomly chosen initial pairs, respectively $O(N/\sqrt{n})$ and $O(nN)$ for $(d)=(2,\ldots,2)$ can be explained as follows. The proof of the theoretical bound starts by bounding
\[
\mathcal{C}_0=\int_0^T\mu(h_t,\zeta_t)\|(\dot h_t,\dot\zeta_t)\|\,dt\leq\sqrt{2}\int_0^T\mu(h_t,\zeta_t)^2\|\dot h_t\|\,dt,
\]
which follows from the fact that $\|\dot\zeta_t\|\leq\mu(h_t,\zeta_t)\|\dot h_t\|$ by the geometric interpretation of the condition number. This last inequality is not sharp in general, and hence one may expect a better behavior of the random linear homotopy method than the one given by the theoretical bound.

\subsection{Equiprobable roots via random homotopy}\label{subsec:equidistribution}
The algorithm constructing a random homotopy has been implemented in two variants:
\begin{enumerate}
\item as described in Section~\ref{sec: random};
\item the initial pair for the linear homotopy is built by taking $(g,e_0)$ in (\ref{equ: conjecture}) and performing a random unitary coordinate transformation (see \cite{Mezzadri} for a stable and efficient algorithm that chooses such a random unitary matrix).
\end{enumerate}

Then the following experiment was conjured to show the equiprobability of the roots at the end of a random homotopy promised by Theorem~\ref{th:complexityhomotopy}: as the target system we take $f=g+\varepsilon h$ where $g$ is as in (\ref{equ: conjecture}), $h$ is chosen randomly in $\S$, and $\varepsilon$ is small. Note that $g$ has a unique non--singular solution which is very well--conditioned, but it also has a whole subspace of degenerate solutions. Hence, $f$ also has a rather well--conditioned solution, and then $\mathcal{D}-1$ isolated, but poorly conditioned ones. One might expect that the random homotopy (2) we have just described (for such a fixed $f$) would be biased to discover the well--conditioned root. Indeed, we obtained numerical evidence that this is not the case: all the solutions seem to be equiprobable.

For $f$ with the degrees $d=(2,2,2)$ and $\varepsilon=0.1$ and several random choices of $g$ we have made experiments with certified tracking procedure making 8000 runs. We experimented with both variants (1) and (2) of choosing the random initial pair. Each experiment resulted in close to 1000 hits for each of 8 roots --- in both variants (1) and (2). This appears to show the conclusion of Theorem \ref{th:complexityhomotopy}, valid for variant (1), and moreover extend it to the case of variant (2).

We can state this experimental result in a more precise way, using Shannon's entropy as suggested in \cite{BeltranPardoFLH}. Assume that we have an algorithm that involves some random choice in its input, and that can produce different outputs $x_1,\ldots,x_l$. Shannon's entropy is by definition the number
\[
H=-\sum_{i=1}^lp_i\log_2(p_i),
\]
where $p_i$ is the probability that the output is $x_i$. It is easy to see that Shannon's entropy of an algorithm is maximal, and equal to $\log_2(l)$, if and only if every output is equally probable. The experimental value of Shannon's entropy for the random algorithm in all experiments described above is in the interval $[2.99,3]$; the maximum, in this case, is $\log_2 8 = 3$.

The exact reason for the modified algorithm (variant (2)) to produce equiprobability of the roots is not understood. This poses a very interesting mathematical question, which together with proving (\ref{eq:SS94B}) would yield a great progress in the understanding of homotopy methods for solving systems of polynomial equations.


\providecommand{\bysame}{\leavevmode\hbox to3em{\hrulefill}\thinspace}
\providecommand{\MR}{\relax\ifhmode\unskip\space\fi MR }
\providecommand{\MRhref}[2]{%
  \href{http://www.ams.org/mathscinet-getitem?mr=#1}{#2}
}
\providecommand{\href}[2]{#2}


\end{document}